\documentclass[a4paper,11pt]{article}

\usepackage{xypic,amsmath,amssymb,latexsym,enumerate}
\usepackage{theorem,eepic,amsbsy}
 \usepackage{geometry,xspace}
\usepackage[curve]{xy}
\usepackage{float}
\usepackage{epsfig}
\usepackage{epic}
\usepackage{overpic}

\usepackage[ps2pdf, bookmarks=true,
            bookmarksnumbered=true, breaklinks=true,
            pdfstartview=FitH, hyperfigures=false,
            plainpages=false, naturalnames=true,
            colorlinks=true,pagebackref=true,
            pdfpagelabels]{hyperref}
\usepackage[usenames]{color}

\newcommand{\quadruple}[4]{\left( \def\arraycolsep{.1cm}
{\begin{array}{ccc} & #3 & \\[-.2cm] #1 & & #4 \\
 [-.2cm] & #2 & \end{array}}
\right)}

\newcommand{\Cat}{\mathsf{Cat}}
\newcommand{\CAT}{\mathsf{CAT}}

\def\END{\mathsf{END}}

\newcommand{\Aut}{\mathsf{Aut}}
\newcommand{\AUT}{\mathsf{AUT}}

\def\K{\mathbb{K}}

\def\2cat{\mathrm{2-}\cat}

\def\Cub{\mathsf{Cub}}

\def\Crs{\mathsf{Crs}}

\def\CRS{\mathsf{CRS}}

\newtheorem{example}{Example}[section]
{\theorembodyfont{\rmfamily}}
{\theorembodyfont{\rmfamily}\newtheorem{Def}[example]{Definition}}
{\theorembodyfont{\rmfamily}}
{\theorembodyfont{\rmfamily}}

\newtheorem{proposition}[example]{Proposition}

\newtheorem{blank}[example]{\hspace{-0.3em}}
{\theorembodyfont{\rmfamily}}

\newenvironment{proof}{\noindent {\bf Proof} }{\hfill $\Box$

\mbox{}}

\newcommand{\directs}[2]{\def\objectstyle{\scriptstyle} \objectmargin={0pt}
\xy
(0,4)*+{}="a",(0,-2)*+{\rule{0em}{1.5ex}#2}="b",(7,4)*+{\;#1}="c"
\ar@{->} "a";"b" \ar @{->}"a";"c" \endxy }

\newcommand{\xdirects}[2]{\def\objectstyle{\scriptstyle} \objectmargin={0pt}
\xy
(0,0)*+{}="a",(0,-6)*+{\rule{0em}{1.5ex}#2}="b",(7,0)*+{\;#1}="c"
\ar@{->} "a";"b" \ar @{->}"a";"c" \endxy }

\newcommand{\sdirects}[2]{\def\objectstyle{\scriptstyle} \objectmargin={0pt}
\xy
(0,2.2)*+{}="a",(0,-2.5)*+{\rule{0em}{1.5ex}#2}="b",(7,2.2)*+{\;#1}="c"
\ar@{->} "a";"b" \ar @{->}"a";"c" \endxy }

\newcommand{\threeaxesb}[3]{\def\objectstyle{\scriptstyle}  \objectmargin={0pt}
\xy
(0,0)*+{}="a",(0,-8)*+{\rule{0em}{1.5ex}#2}="b",(10,0)*+{\;#1}="c",
(8,7)*+{\;#3}="d" \ar@{->} "a";"b" \ar @{->}"a";"c"  \ar
@{->}"a";"d"\endxy }

\newcommand{\sq}{\mbox{\rule{0.08em}{1.7ex}\hspace{-0.00em}\rule{0.7em}{0.2ex}\hspace{-0.7em}\rule[1.54ex]{0.7em}{0.2ex}\hspace{-0.03em}\rule{0.08em}{1.7ex}}}

\newcommand{\io}{^{-1}}

\def\<{\langle \! \langle}
\def\>{\rangle \! \rangle}

\newcommand{\Ob}{\mathrm{Ob}}

\def\Ob{\operatorname{Ob}}

\newcommand{\C}{\mathsf{C}}

\def\A{\mathsf{A}}

\def\leq{\leqslant}

\newcommand{\labto}[1]{\stackrel{#1}{\longrightarrow}}

\begin{document}
\title{Possible connections between\\
whiskered categories and groupoids, many object Leibniz algebras, \\
automorphism structures and local-to-global questions}
\author{Ronald Brown\thanks{School of Computer Science, University of Wales, Bangor, Gwynedd LL57 1UT, U.K., r.brown
 `at ' bangor.ac.uk, www.bangor.ac.uk/r.brown}}
\maketitle
\begin{abstract}
We define the notion of whiskered categories and groupoids, showing
that whiskered groupoids have a commutator theory. So also do
whiskered $R$-categories, thus answering questions of what might be
`commutative versions' of these theories. We relate these ideas to
the theory of Leibniz  algebras, but the commutator theory here does
not satisfy the Leibniz  identity. We also discuss potential
applications and extensions, for example to resolutions of monoids.
\end{abstract}

\section*{Introduction}
The notion of commutativity is standard for monoids and groups, but
seems to be lacking for categories and groupoids. Similarly, there
is a notion of Lie bracket $[a,b]=ab-ba$ for elements $a,b$ of an
associative algebra $A$, but this seems to have no parallel for the
case of additive categories, which can be seen as algebras with many
objects.

The aim of this paper to introduce the extra structure of whiskering
in these situations so that we can discuss commutativity,
commutators, `Lie brackets',  and related questions.

It was originally expected that a generalised Lie  bracket would
satisfy axioms analogous to the Jacobi, or later the Leibniz ,
identity. However it turns out that in these whiskered situations
the rules satisfied by the commutator or bracket which occur  are
best described in a cubical background, and so  Section
\ref{sec:cubes} is devoted to this account. More exploitation of
cubes in a category is given in the first section of \cite{BL2}.

Section \ref{sec:whiskering} gives the definition of a whiskered
category. Section \ref{sec:commutators} introduces commutators in a
whiskered groupoid and examines their basic properties. Section
\ref{sec:Rcat} discusses the properties of `Lie brackets' in a
whiskered $R$-category. Section \ref{sec:braidedcrossedcompl} shows
how whiskered categories and groupoids arise from braidings on
crossed complexes.  Section \ref{sec:resolmonoids} discusses
possibilities for resolutions of monoids. Section
\ref{whiskerandcubical} discusses the role of the cubical theory in
this area.

\section{Squares and cubes in categories}\label{sec:cubes}
Let $\C$ be a category, with set of objects written $\C_0$. The set
of arrows of $\C$ is written $\C_1$. We will write the composition
of $a: x \to y$ and $b: y \to z$ in the algebraic fashion as $ab: x
\to z$ or $a \circ b$, since this is more convenient for the
algebraic work here.

\begin{Def}
We write $I$ for the ordered set $\{-,+\}$ with $- < +$, also
regarded as a category. A {\it square, or $2$-cube,} in the category
$\C$ is a functor $f: I^2 \to \C$ and this is written as a diagram
\begin{equation} \label{diag:sq}
\def\labelstyle{\textstyle} \vcenter{\xymatrix@=3pc{ x \ar [r] ^{\partial^-_1 f} \ar
[d] _ {\partial^-_2 f} &  \ar [d]^{\partial^+_2 f}
\\  \ar [r] _{\partial^+_1 f} & y}}\quad \directs{2}{1}
\end{equation}
We define $sf=x , tf=y$ as in the diagram. The squares in $\C$ form
with the obvious compositions a double category $\sq \C$, with
compositions $\circ_1. \circ_2$.  \end{Def}

\begin{Def}
If further $\C$ is a groupoid we define
$$\delta f= (\partial ^+_2 f)\io (\partial ^-_1 f)\io (\partial ^-_2 f)(\partial ^+_1f ), $$
which clearly belongs to $\C_1(y,y)$.  \end{Def} We now turn to the
additive case.
\begin{Def}
If $\C$ is an additive category and $f: I^2 \to \C$ is a square in
$\C$ we define
$$\Delta f= -(\partial ^-_2 f)(\partial ^+_1 f)+ (\partial^ -_1 f)(\partial
^+_2 f)$$ which clearly belongs to $\C(sf,tf)$.  \end{Def} If $\C$
is additive then $\sq \C$ obtains additional partial additive
structures as follows:
\begin{alignat*}{2}
  \quadruple{a}{b}{c}{d} &+_1  \quadruple{a}{b'}{c'}{d}&&=
\quadruple{a}{b+b'}{c+c'}{d}\\
 \quadruple{a}{b}{c}{d} &+_2  \quadruple{a'}{b}{c}{d'}&&=
\quadruple{a+a'}{b}{c}{d+d'}. \tag*{$\Box$}
\end{alignat*}
Note that
$$\Delta(\alpha +_i \beta)= \Delta(\alpha)+ \Delta(\beta) $$
for $i=1,2$ and assuming $\alpha+_i \beta$ is defined.

We record the following for later use.
\begin{proposition}
(i) Let $\C$ be a groupoid and $\sq \C$ its double category of
squares. Then for $\alpha, \beta, \gamma \in \sq \C$ such that
$\alpha \circ_1\beta, \alpha \circ_2 \gamma$ are defined:
\begin{align*}
\delta(\alpha \circ _1 \beta) & = (\delta \beta) (\delta \alpha)^
{\partial^+_2 \beta},\\
\delta(\alpha \circ _2 \gamma )& = (\delta \alpha)^ {\partial^+_1
\gamma} (\delta \gamma).
\end{align*}
(ii) Let $\C$ be an additive  category and $\sq \C$ its double
additive category of squares. Then for $\alpha, \beta, \gamma \in
\sq \C$ such that $\alpha \circ_1\beta, \alpha \circ_2 \gamma$ are
defined:
\begin{align*}
\Delta(\alpha \circ _1 \beta) & = (\Delta \alpha)
(\partial^+_2 \beta)+ (\partial^-_2) \alpha)(\Delta \beta) ,\\
\Delta(\alpha \circ _2 \gamma )& = (\Delta \alpha)(\partial^+_2
\gamma) +  (\partial^-_1 \gamma)(\Delta \gamma).
\end{align*}
\end{proposition}
The proofs are straightforward. The formulae (i) in the above
Proposition are related to formulae occurring in the relations
between double groupoids and crossed modules, see for example
\cite[Section 6.6]{BHS}.

A {\it $3$-cube} in the category $\C$ is a functor $f:I^3 \to \C$.
\begin{proposition} \label{prop:morphfordDelta}
 If $f$ is a $3$-cube in the groupoid $\C$ then we have the rule
\begin{equation} \label{comm:cubdir0}
 (\delta \partial ^-_3 f)^{u_3} (\delta \partial ^+_2 f) (\delta
\partial^-_1 f)^{u_1}= (\delta \partial ^+_1 f ) (\delta \partial ^-_2 f)^{u_2}
(\delta \partial ^+_3 f),
\end{equation}
where $u_1= \partial^+_2\partial^+ _3f, u_2= \partial^+_1\partial^+
_3f, u_3= \partial^+_1\partial^+ _2f$ and $a^b=b \io ab$.
\end{proposition}
\begin{proof}
It is convenient to label the edges of the cube as follows:
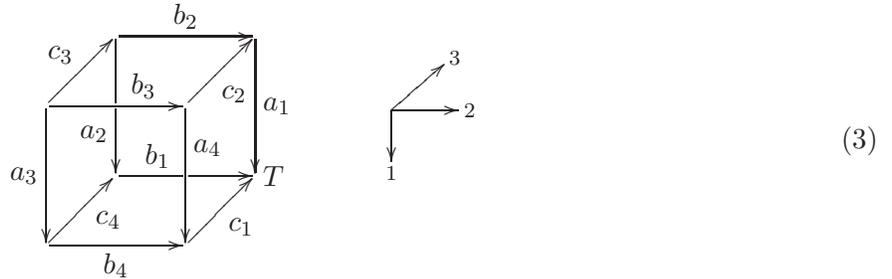
\begin{equation}\label{comm:cubdir-abc}\def\labelstyle{\textstyle}
{{\objectmargin{0.1pc} \diagram &  \rrto^{b_2} \xline '[1,0]
[2,0]|>\tip \ar @{} [dd]_(0.7){a_2}&&\ddto^{a_1} ^(1.0)T\\
\urto^{c_3} \rrto^(0.7){b_3} \ddto_{a_3} && \urto_(0.4){c_2}
\ddto^(0.3){a_4}&
\\&\ar @{}[rr]^(0.3){b_1}\xline'[0,1] [0,2]|>\tip &&  \\ \rrto_{b_4} \urto_(0.6){c_4} &&
\urto_{c_1} & \enddiagram }}\qquad \raisebox{3ex}{
\threeaxesb{2}{1}{3}}\end{equation} so that $a_1=u_1, b_1=u_2,
c_1=u_3$. Then both sides of equation \eqref{comm:cubdir0} reduce
to:
\begin{equation*} a_1\io b_2\io c_3\io a_3b_4c_1.
\end{equation*}
\end{proof}

\begin{proposition} \label{prop:Deltaassoc}
 Let $f$ be a $3$-cube in an additive category $\C$. Let $$a_3= \partial^-_2 \partial^-_3
f, a_1=\partial^+_2 \partial^+_3 f,   b_3=\partial^-_1 \partial^-_3,
b_1= \partial^+_1 \partial^+_3 f, c_3=\partial^-_1
\partial^-_2 f,c_1= \partial^+_1 \partial^+_2 f .$$ Then
$$\Delta\!\quadruple{a_3}{\Delta \partial^+_1 f}{\Delta \partial^-_1
f}{a_1}= \Delta\!\quadruple{b_3}{\Delta \partial^+_2 f}{\Delta
\partial^-_2 f}{b_1}-\Delta\!\quadruple{c_3}{\Delta \partial^+_3 f}{\Delta \partial^-_3
f}{c_1}.$$
\end{proposition}
\begin{proof}
We label the edges of the cube $f$ as in diagram
\eqref{comm:cubdir-abc}. Then the definitions imply that
\begin{align*}
\Delta\!\quadruple{a_3}{\Delta \partial^+_1 f}{\Delta \partial^-_1
f}{a_1}&=
\Delta\!\quadruple{a_3}{-b_4c_1+c_4b_1}{-b_3c_2+c_3b_2}{a_1}\\
&= a_3b_4c_1-a_3c_4b_1-b_3c_2a_1+c_3b_2a_1, \\
\Delta\!\quadruple{b_3}{\Delta \partial^+_2 f}{\Delta
\partial^-_2 f}{b_1}&= b_3a_4c_1-b_3c_2a_1-a_3c_4b_1+c_3a_4b_1, \tag*{similarly,}\\
\Delta\!\quadruple{c_3}{\Delta \partial^+_3 f}{\Delta \partial^-_3
f}{c_1}&= c_3a_2b_1-c_3b_2a_1-a_3b_4c_1+b_3a_4c_1,
\end{align*}
from which the result follows.
\end{proof}

\section{Whiskering}\label{sec:whiskering}
\begin{Def}
 A {\it whiskering} on a category $\C$ consists of operations
\begin{align*}
  m_{ij}: \C_i \times \C_j \to \C_{i+j}, \quad i,j=0,1,\;  i+j \leq 1,
\end{align*}
satisfying the following axioms:
\begin{enumerate}[\noindent {Whisk} 1.)]
\item $m_{00}$ gives a monoid structure on $\C_0$ with identity written $1$ and multiplication written as juxtaposition;
 \item
$m_{01},m_{10}$ give respectively left and right actions of the
monoid $\C_0$ on the category $\C$, in the sense that: \item if $x
\in \C_0$ and $ a: u \to v$ in $\C_1$, then $x.a: xu \to xv$ in
$\C$, so that $$1.a=a, (xy).a=x.(y.a), $$$$ x.(a\circ b)=(x.a)\circ
(x.b), x.1_y=1_{xy};$$\item analogous rules hold for the right
action; \item $x.(a.y)=(x.a).y$,
\end{enumerate}
for all $x,y,u,v \in \C_0, a,b \in \C_1$.

A {\it whiskered category} is a category with a whiskering.
\end{Def}

\begin{Def}
Let $\C$ be a category. A {\it bimorphism} $m:(\C,\C) \to \sq \C$
assigns to each pair of morphisms $a,b \in \C$ a square $m(a,b) \in
\sq \C$  such that if $ad, bc $ are  defined in $\C$ then
\begin{align*}
m(ad,c) &= m(a,c) \circ_1 m(d,c),\\
m(a,bc)&= m(a,b) \circ_2 m (a,c). \tag*{$\Box$}
\end{align*}
\end{Def}
The following is easy to prove.
\begin{proposition}
 If $\C$ is a whiskered category then a bimorphism
$$*: (\C,\C) \to \sq \C$$
is defined for $a: x \to y, b: u \to v $ by
$$a*b= \quadruple{a.u}{y.b}{x.b}{a.v}. $$
\end{proposition}

 If $\C$ is a whiskered category, then two multiplications $l,r$
on the set $\C_1$ may be defined by:\\ if $a:x\to y, b:u \to v$,
then
$$ l(a,b):=(a.u)\circ (y.b),\qquad  r(a,b):= (x.b) \circ (a.v), \quad  $$
as shown in the diagram \begin{equation} \label{eq:noncommsquare}
\def\labelstyle{\textstyle} \vcenter{\xymatrix@=3pc{ xu \ar [r] ^{x.b} \ar
[d] _ {a.u} & xv \ar [d]^{a.v}
\\ yu \ar [r] _{y.b} & yv}}\quad \directs{2}{1}\end{equation}
\begin{proposition}
If $l(a,b)=r(a,b)$ for all $a,b \in \C_1$, then the multiplication
$(a,b)\mapsto a.b$  given by this common value makes $\C$ into a
strict monoidal category.
\end{proposition}
\begin{proof}
Suppose also $c:y \to z,d:v \to w$. Then the commutativity of the
diagram \begin{equation} \label{diag:interch}
\def\labelstyle{\textstyle} \vcenter{\xymatrix@=3pc
 { xu \ar [d] _{a.u} \ar [r] ^ {x.b} & xv \ar [d]|{a.v}\ar [r] ^{x.d} &xw \ar
 [d] ^{a.w} \\ yu \ar [r] |{y.b} \ar [d] _{c.u}& yv \ar [r] |{y.d}\ar [d] |{c.v} & yw \ar [d]
 ^ {c.w}\\zu \ar [r] _{z.b} &zv \ar [r] _{z.d} & zw}}\end{equation}
yields a verification of the interchange law for $.$ and $\circ$.

The verifications of the laws  for associativity and the identity
are trivial.
\end{proof}

In the case given by this proposition we say $(\C,m)$ is a {\it
commutative} whiskered category.

In general though the interchange law is not satisfied and what we
have is a {\it sesquicategory} as considered in
\cite{Stell,tonkshess:model-2cat,HeJohnson,street-catstr}. It is
notable that a majority of writing on weak or lax forms of
$n$-categories, see for example \cite{Cheng-Gurski} and the
references there,  is in the globular format and assumes a strict
interchange law. However, as we discuss in section
\ref{sec:cubical}, there is a case for a cubical approach and
failure of the interchange law is interesting, is potentially
controllable with some special structures, and seems to arise in
geometric situations. Indeed, since one overall aim of higher
category theory is to model homotopy types in a useful way, the fact
that weak, pointed homotopy types are modelled by simplicial groups
should take account of the complex structures this entails, see for
example \cite{CC}. Such truncated homotopy types are also modelled
by the strict cat$^n$-groups, \cite{Lod82}, and this allows for some
calculation, \cite{BL87}, and extension of classical homotopical
results such as the Blakers-Massey and Barratt-Whitehead theorems
connectivity and critical group theorems, and the relative Hurewicz
theorem, \cite{BL2}. Relations between such strict and some weak
models are considered in \cite{Paoli-n-types}.

%\newpage

\section{Commutators in whiskered groupoids}\label{sec:commutators}
In the case $\C$ is a whiskered groupoid, and $a,b \in \C_1$ as in
the proposition, we define the {\it commutator}
\begin{equation}\label{eq:gpdcommutator}
[a,b] = \delta(a*b) .\end{equation} Notice that this definition
requires a convention as to starting point and direction around the
square.

 It is interesting to see how the usual laws for commutators
generalise to these commutators in a whiskered groupoid. We now
abbreviate $a\circ b$ to $ab$.

One of the easiest rules for the usual commutators fails in this
context.  Thus when $a: x \to y$  we find
\begin{equation}
  [a,a]= (a.y)\io(x.a)\io(a.x)(y.a),
\end{equation}
so that in general $[a,a] \ne 1$. Similarly $[a,b] \ne [b,a]\io$.

\begin{proposition}
The whiskered groupoid $\C$ satisfies the rules $[a,a]=1,
[a,b]=[b,a]\io $ for all $a,b \in \C$ if and only if the monoid
$\C_0$ is commutative and the action satisfies $x.a=a.x$ for all $a
\in \C, x \in \C_0$.
\end{proposition}

A biderivation rule for commutators carries over to this situation:
\begin{proposition} Let $a:x \to y, b: u\to v , c: y \to z, d: v \to w$ in $\C_1$. Then
\begin{align*}
[ac,b]& = [a,c]^{c.v} \; [c,b],\\
[a,bd]&= [a,d][a,b]^{y.d}.
\end{align*}
\end{proposition}
The proof is straightforward, by referring to Proposition
\ref{prop:morphfordDelta} or diagram \eqref{diag:interch}.

This result suggests the possibility of a nonabelian tensor product,
see \cite{BL87}.

In the case of groups there is a well known law on commutators which
with our convention  reads as \begin{equation} \label{comm:rule1}
 [a,b]^c\; [a,c][b,c]^a= [b,c][a,c]^b\; [a,b],
\end{equation} where $[a,b]=a\io b\io ab, \; x^a=a \io x a$.
This equation may be viewed cubically as:
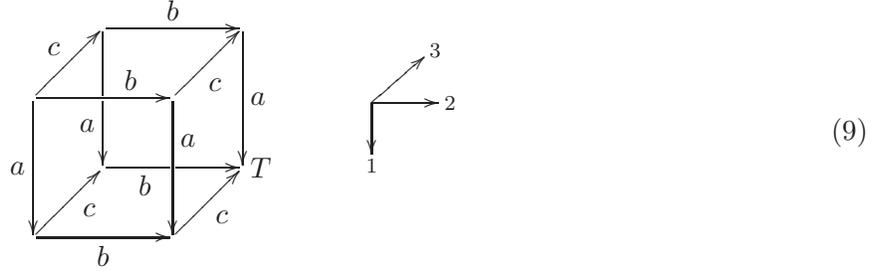
\begin{equation}\label{comm:cubdir1}\def\labelstyle{\textstyle}
{{\objectmargin{0.1pc} \diagram &  \rrto^b \xline '[1,0] [2,0]|>\tip
\ar @{} [dd]_(0.7)a&&\ddto^a ^(1.0)T\\ \urto^c \rrto^(0.7)b \ddto_a
&& \urto_(0.4)c \ddto^(0.3)a &
\\&\ar @{}[rr]_(0.3)b\xline'[0,1] [0,2]|>\tip &&  \\ \rrto_b \urto_(0.6)c &&  \urto_c
& \enddiagram }}\qquad \raisebox{3ex}{
\threeaxesb{2}{1}{3}}\end{equation} With the directions shown by the
axes the rule \eqref{comm:rule1} can be read as:
\begin{equation}\label{comm:rule2}
(\partial^-_3)^c \; (\partial^+_2)(\partial^-_1)^a=
(\partial^+_1)(\partial^-_2)^b \; (\partial^+_3).\end{equation} Thus
we see the role of the exponents is to bring the appropriate face to
have base point at the far right hand corner labelled $T$.

To obtain an analogous identity in our whiskered situation we again
need to write out a cubical model. Thus if $a:x\to y, b: u \to v, c:
z \to w$ then $a*b*c$ may be seen as the 3-cube in $\C$ given as
follows:\begin{equation} \label{cubediag:xuz}
\def\labelstyle{\textstyle}{{\objectmargin{0.1pc} 
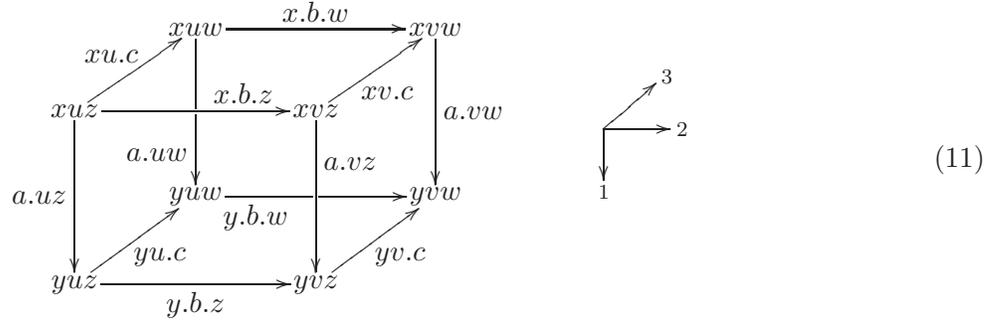
\diagram &  xuw \ar@{}[dd]_(0.75){a.uw}\rrto^{x.b.w} \xline '[1,0]
[2,0]|>\tip &&xvw\ddto^{a.vw}\\xuz \urto^{xu.c} \rrto^(0.7){x.b.z}
\ddto_{a.uz} && \urto_(0.4){xv.c} \ddto^(0.3){a.vz} xvz&
\\&yuw \ar@{}[rr]_(0.25){y.b.w}\xline'[0,1] [0,2]|>\tip && yvw \\ yuz\rrto_{y.b.z} \urto_{yu.c} && yvz \urto_{yv.c}
& \enddiagram }}\qquad \raisebox{3ex}{
\threeaxesb{2}{1}{3}}\end{equation}

We then have the following description of the faces.
\begin{proposition}\label{lem:facesofcube}
The faces of the cube then determine  commutators as follows:
\begin{alignat*}{2}
\partial^-_1&= x.[b,c], \quad& \partial^+_1 &= y.[b,c],\\
\partial^-_2&= [a,u.c], & \partial^+_2 &= [a,v.c],\\
\partial^-_3&= [a,b].z, & \partial^+_3 &= [a,b].w.
\end{alignat*}
\end{proposition}
\begin{proposition}
These commutators satisfy the rule
\begin{equation*}\label{comm:rule3}
 ([a,b].z)^{yv.c}\;
([a,v.c])(x.[b,c])^{a.vw}= (y.[b,c])([a,u.c])^{y.b.w} \;
([a,b].w).\end{equation*}
\end{proposition}
\begin{proof}
This follows immediately from Proposition \ref{prop:morphfordDelta}.
\end{proof}

\section{Whiskered $R$-categories and Leibniz
algebras}\label{sec:Rcat} Let $R$ be a commutative ring and let $\C$
be  a whiskered category. Then we can form the category $R[\C]$ with
the same objects as $\C$ but with $R[\C](x,y)$ the free $R$-module
on $\C(x,y)$ for all $x,y \in \C_0$. The the action of $\C_0$ on
$\C$ extends to an action of $\C_0$ on $R[\C]$ which is bilinear in
the sense that
\begin{equation}
  x.(ra + a') = r(x.a) +a',\quad  (ra+a') .x = r(a.x) + a'.x,
\end{equation}
for all $r \in R, x \in C_0, $ and $a,a'$ in $\C$ with the same
source and target. In such case we say $R[\C]$ is a {\it whiskered
$R$-category}.

Suppose now that $\A$ is a whiskered $R$-category. We can
analogously to the above define the Lie bracket of $a: x \to y$ and
$b: u \to v$ by
\begin{equation}\label{Leibniz comm} [a,b] = \Delta(a*b)= -((a.u) (y.b) )+
(x.b)(a.v).\end{equation} Again we see that
$$[a,a] = -((a.x)(y.a)) + ((x.a)  (a.y))$$
so that in general $[a,a] \ne 0$. This suggests that we might have
not a Lie algebra but a Leibniz  algebras, \cite{loday-Leibniz },
which in the usual case requires the {\it Leibniz  identity}
\begin{equation}
[[a,b],c]= [a,[b,c]] + [[a,c],b].
\end{equation}
However instead we have the equation as follows.
\begin{proposition}
 If $\C$ is a whiskered $R$-category and $a:x \to y, b: u: \to v, c:
z \to w$ in $\C_1$, then
\begin{equation*}
[[a,b],c]- [a,[b,c]] =
\Delta\quadruple{[a,u.c]}{y.b.w}{x.b.z}{[a,v.c]}.
\end{equation*}
\end{proposition}
\begin{proof}
This follows immediately from Propositions \ref{prop:Deltaassoc} and
 \ref{lem:facesofcube}, as the latter description of the faces of
the cube holds also for Lie brackets as well as for commutators.
\end{proof}

Thus we have not found a solution to the problem raised by Loday in
\cite{loday-Leibniz } of the existence of what he calls a
`coquecigrue', i.e. a group like structure whose representations
form a Leibniz  algebra. That paper, and others, are also interested
in the smooth case, and are asking for a differentiable coquecigrue
which has an associated Leibniz  algebra, analogous to the way a Lie
group has an associated Lie algebra.

This also leaves open the question of properties such as the
Poincar\'e-Birkhoff-Witt theorem, analogous to the classical case as
discussed for example in \cite{Higgins-baer-lie}. A relevant paper
is \cite{loday-pir-envoloping}.

A further question is whether these ideas are useful for
generalising  the theory of Lie algebroids of Lie groupoids as in
\cite{Mackenzie2005} to the case of Lie 2-groupoids and other Lie
multiple groupoids.

\section{Braided crossed complexes and
automorphisms}\label{sec:braidedcrossedcompl} The category $\Cat$ of
small categories is cartesian closed with an exponential law of the
form
$$\Cat(A \times B, C) \cong \Cat(A, \CAT(B,C))$$
for all small categories $A,B,C$. It follows that for any small
category $C$, $\END(C)= \CAT(C,C)$ is a monoid in $\Cat$ which has a
maximal subgroup object which we call $\AUT(C)$, the {\it actor} or
{\it automorphism object} of $C$. An exposition of these matters is
in \cite{BHS}.

The above facts have analogues in any cartesian closed category.

In the case $C$ is a groupoid, then $\AUT(C)$ is equivalent to a
crossed module $\xi: S(C) \to \Aut(C)$ where $S(C)$ is the set of
bisections of $C$, i.e. sections $\sigma$ of the source map $s$ such
that $t\sigma$ is a bijection on $\Ob(C)$. The set $S(C)$ has a
group structure with the Ehresmannian composition $\tau \circ \sigma
(x) = \tau(t\sigma x) \sigma x$, for $x \in \Ob(C)$. The
automorphisms in the image of $\xi$ are called {\it inner
automorphisms} of $C$.

The situation is rather different in a monoidal closed category. For
example, Brown and Gilbert in \cite{Brown-Gilbert} considered the
monoidal closed category of crossed modules of groupoids. This was
deduced from the  monoidal closed structure on the category $\Crs$
of crossed complexes, with an exponential law of the form
\begin{equation}
\Crs(A \otimes B, C) \cong \Crs(A, \CRS(B,C)).
\end{equation}
constructed by Brown and Higgins in \cite{BH87}.  The monoidal
closed structure is used in an essential way to formulate a model
structure for the homotopy of crossed complexes as in \cite{BG89}.
This exponential law implies that $\END(C)= \CRS(C,C)$ is a monoid
in $\Crs$ with respect to $\otimes$, but the concept of group object
with respect to $\otimes$ is not meaningful, so we have to take a
different approach to obtain a candidate for the actor of a crossed
complex.

Now $\CRS(A,B)_0= \Crs(A,B)$. So we define $\AUT(C)$ to be the full
subcrossed complex of $\END(C)$ on the set $\Aut(C) \subseteq
\END(C)_0$. Clearly $\AUT(C)$ is a monoid object in $\Crs$ with
respect to $\otimes$.

A generalisation of a construction in \cite{Brown-Gilbert} is to
form the simplicial nerve $N^\Delta(\AUT(C))$. Recall that for a
crossed complex $A$ the simplicial nerve $N^\Delta(A)$ is defined to
be in dimension $n$
\begin{equation}
N^\Delta(A)_n= \Crs(\Pi \Delta^n, A),
\end{equation}
where $\Delta^n$ is the $n$-simplex and $\Pi$ gives the fundamental
crossed complex. The crossed complex $\Pi \Delta^n$ is constructed
in \cite{BS:normalisation} directly from the monoidal structure on
$\Crs$.

Now there is a crossed complex Eilenberg-Zilber theorem due to Tonks
in \cite{tonks:EZ}. This gives an Alexander-Whitney type diagonal
map \begin{equation} AW: \Pi \Delta^n \to (\Pi \Delta^n) \otimes
(\Pi \Delta^n).
\end{equation}
So given a morphism $m:A \otimes A \to A$ we get a `convolution'
product $f*g\in N^\Delta(A)_n$ of $f,g\in N^\Delta(A)_n$ as the
composition
\begin{equation}
 \Pi \Delta^n \labto{AW}(\Pi \Delta^n) \otimes
(\Pi \Delta^n) \labto{f \otimes g} A \otimes A \labto{m} A.
\end{equation}
The  properties of the map $AW$ as given in \cite{tonks:EZ} and the
properties of $m$, including the fact than $A_0$ is a group, then
imply that $N^\Delta \AUT(C)$ has an induced structure of simplicial
group. This is the justification for considering $\AUT(C)$ as a
candidate for the actor (automorphism structure) of a crossed
complex.

In the terminology of \cite{Brown-Gilbert}, we would call $\AUT(C)$
a braided, regular crossed complex. See also \cite{arvasi-ulualan}
for related material in the crossed module case. In
\cite{Baues-Tonks-twisted}, crossed differential algebras are called
crossed chain algebras.

The theory of \cite{Brown-Gilbert} was applied to the case of
crossed modules of groups and the corresponding application to
crossed modules of groupoids was worked out in
\cite{BrIc-automorph}.  By working entirely in these crossed modules
of groupoids, some of the proofs seem detailed and unintuitive, and
we felt that they would be better in terms of 2-groupoids. However,
homotopies of 2-groupoids are more complicated than homotopies of
crossed modules, and this project was not completed. Related work
considering automorphisms of 2-groups using 2-groupoids is in
\cite{Roberts-Schreiber}. Crossed complexes are equivalent to
globular $\infty$-groupoids (sometimes called $\omega$-groupoids) by
work of \cite{BH81:inf}.

Whitehead's work on operators on relative homotopy groups in
\cite{Whitehead:operators} was continued by Hu in
\cite{Hu-automorphisms}; it may be worth taking another look at
these matters from a modern and broader perspective.

The work as given in \cite{BrIc-automorph} was necessary for the
work on 2-dimensional holonomy in \cite{Br-Ic-2holonomy}.
Originally, we naively conjectured that as a foliation gave rise to
a holonomy groupoid, so a double foliation would give rise to a
holonomy double groupoid. This was not achieved in
\cite{Br-Ic-2holonomy}, and instead we obtained only (or so it
seemed) a 2-crossed module from situations in this area, following
ideas in \cite{Aof-Br}, but generalising local sections to
homotopies. Perhaps this is inevitable, and local interchange laws
do not necessarily lead to global interchange laws, because of the
influence of non local features, just as a bundle can be locally
trivial but not trivial. The analysis of this distinction needs an
appropriate structure, which in the case of bundles has been known
since the work of Ehresmann.  Thus  the control of this lack of
global interchange law given by, say,  a 2-crossed module could be
important. The notion of locally topological 2-crossed module (of
groupoids?) has not yet been considered!

An alternative to the simplicial theory indicated above, and which
has not been worked on in this context, is to consider the cubical
nerve $N^\Box A$ of a crossed complex, with values in the category
of cubical sets. This construction gives the equivalence between the
category $\Crs$ and that of cubical $\omega$-groupoids with
connections, see \cite{BH81:algcub}. Thus we should perhaps consider
the automorphism theory not in crossed complexes but in the natural
home for homotopies and tensor products,  namely the monoidal closed
category of cubical $\omega$-groupoids with connections, see
\cite{BH87}.

\section{Resolutions of monoids?}\label{sec:resolmonoids}
The use of crossed differential algebras suggests a possibility for
resolutions of monoids. We know that a quotient of a monoid is
described by a congruence, which is an equivalence relation in the
category of monoids.  For obtaining free objects it is natural
therefore to consider groupoid objects in the category of monoids.

There is a different way of considering this question. Let $M$ be a
monoid. Then $M$ defines a crossed differential algebra $A=\K(M,0)$
which is $M$ in dimension 0, trivial otherwise and with monoid
structure with respect to $\otimes$ which of course is trivial
except in dimension 0, given by the multiplication on $M$.

Of course $A$ is not free, or cofibrant, in any sense. The category
of crossed differential algebras has a homotopy structure, as shown
by Tonks in \cite{tonks:I-structure}; the paper
\cite{riehl-weakfactoris} is also relevant. So it is interesting to
replace $A$ by a cofibrant object up to weak equivalence. I have not
done the work on this, but possibilities are as follows.

First one chooses a set $X$ of generators of $M$ as a monoid, and
forms the free monoid $X^*$ on $X$ and its associated map $f: X^*
\to M$. The next step is presumably related to work of Porter in
\cite{Porter-gummersbach},  Heyworth and Johnson in \cite{HeJohnson}
and possibly to that of \cite{tonkshess:model-2cat}. It seems one
should choose the free whiskered groupoid on generators of the
congruence given by $f$.

More globally, the methods of
\cite{tonks:I-structure,riehl-weakfactoris} suggest that there is a
cofibrant object in the category of crossed differential algebras
extending $X^*$ and with a weak equivalence to $A$.

Related notions are also in \cite{Gilbert-monoidpres}.

\section{Whiskered categories and cubical
theory}\label{whiskerandcubical} \label{sec:cubical} The context of
crossed complexes is a candidate for the groupoid theory but not for
the category case. There is an argument for the monoidal closed
category, say $\Cub\Cat$, of cubical $\omega$-categories with
connections defined and developed in \cite{ABS}. In this category a
monoid object $A$ with respect to $\otimes$ has as its truncation
$tr^1 A$ exactly a whiskered category. But $tr^2 A$ also contains
the not necessarily commutative square given in the diagram
\eqref{eq:noncommsquare}, and a filler of it, namely the product
$ab\in A_2$ in the monoid structure.

Of course the main result of \cite{ABS} is the equivalence of
$\Cub\Cat$ with the more usual category of globular
$\omega$-categories. The advantage of the cubical case is as usual
the easy definitions of multiple compositions, and of tensor
products, and these are the basis of the topological applications of
the cubical higher homotopy groupoid of a filtered space, see the
survey \cite{Brown-handbook}.

\end{document}